\documentclass[12pt]{article}

\usepackage{amsmath,amssymb,amsbsy,amsfonts,amsthm,latexsym,
                            amsopn,amstext,amsxtra,euscript,amscd}

\begin{document}

\newtheorem{theorem}{Theorem}
\newtheorem{lemma}[theorem]{Lemma}
\newtheorem{claim}[theorem]{Claim}
\newtheorem{cor}[theorem]{Corollary}
\newtheorem{prop}[theorem]{Proposition}
\newtheorem{definition}{Definition}
\newtheorem{question}[theorem]{Open Question}

%%%%%%%%%%%%%%%%%%%%%%%%%
% Alphabet calligraphic %
%%%%%%%%%%%%%%%%%%%%%%%%%
\def\cA{{\mathcal A}}
\def\cB{{\mathcal B}}
\def\cC{{\mathcal C}}
\def\cD{{\mathcal D}}
\def\cE{{\mathcal E}}
\def\cF{{\mathcal F}}
\def\cG{{\mathcal G}}
\def\cH{{\mathcal H}}
\def\cI{{\mathcal I}}
\def\cJ{{\mathcal J}}
\def\cK{{\mathcal K}}
\def\cL{{\mathcal L}}
\def\cM{{\mathcal M}}
\def\cN{{\mathcal N}}
\def\cO{{\mathcal O}}
\def\cP{{\mathcal P}}
\def\cQ{{\mathcal Q}}
\def\cR{{\mathcal R}}
\def\cS{{\mathcal S}}
\def\cT{{\mathcal T}}
\def\cU{{\mathcal U}}
\def\cV{{\mathcal V}}
\def\cW{{\mathcal W}}
\def\cX{{\mathcal X}}
\def\cY{{\mathcal Y}}
\def\cZ{{\mathcal Z}}

%%%%%%%%%%%%%%%%%%%%%%%
% Alphabet blackboard %
%%%%%%%%%%%%%%%%%%%%%%%
\def\A{{\mathbb A}}
\def\B{{\mathbb B}}
\def\C{{\mathbb C}}
\def\D{{\mathbb D}}
\def\E{{\mathbb E}}
\def\F{{\mathbb F}}
\def\G{{\mathbb G}}
\def\H{{\mathbb H}}
\def\I{{\mathbb I}}
\def\J{{\mathbb J}}
\def\K{{\mathbb K}}
\def\L{{\mathbb L}}
\def\M{{\mathbb M}}
\def\N{{\mathbb N}}
\def\O{{\mathbb O}}
\def\P{{\mathbb P}}
\def\Q{{\mathbb Q}}
\def\R{{\mathbb R}}
\def\S{{\mathbb S}}
\def\T{{\mathbb T}}
\def\U{{\mathbb U}}
\def\V{{\mathbb V}}
\def\W{{\mathbb W}}
\def\X{{\mathbb X}}
\def\Y{{\mathbb Y}}
\def\Z{{\mathbb Z}}

\def\E{{\mathbf E}}
\def\Fp{\F_p}
\def\ep{{\mathbf{e}}_p}

\def\scr{\scriptstyle}
\def\\{\cr}
\def\({\left(}
\def\){\right)}
\def\[{\left[}
\def\]{\right]}
\def\<{\langle}
\def\>{\rangle}
\def\fl#1{\left\lfloor#1\right\rfloor}
\def\rf#1{\left\lceil#1\right\rceil}
\def\le{\leqslant}
\def\ge{\geqslant}
\def\eps{\varepsilon}
\def\mand{\qquad\mbox{and}\qquad}

\newcommand{\comm}[1]{\marginpar{%
\vskip-\baselineskip %raise the marginpar a bit
\raggedright\footnotesize
\itshape\hrule\smallskip#1\par\smallskip\hrule}}

\def\xxx{\vskip5pt\hrule\vskip5pt}

%%%%%%%%%%%%%%%%%%
%% PAPER BEGINS %%
%%%%%%%%%%%%%%%%%%

\title{\bf Sato--Tate, cyclicity, and divisibility statistics
on average for elliptic curves of small height}

\author{
{\sc William D.\ Banks} \\
{Department of Mathematics, University of Missouri} \\
{Columbia, MO 65211 USA} \\
{\tt bbanks@math.missouri.edu} \\
\and
\\
{\sc Igor E. Shparlinski} \\
{Department of Computing, Macquarie University} \\
{Sydney, NSW 2109, Australia} \\
{\tt igor@ics.mq.edu.au}}

\date{\today}
\pagenumbering{arabic}

\maketitle

\begin{abstract}
We obtain  asymptotic formulae for the number of primes $p\le x$
for which the reduction modulo $p$ of the elliptic curve
$$
\E_{a,b}~:~Y^2 = X^3 + aX + b
$$
satisfies certain ``natural'' properties, on average  over
integers $a$ and $b$ such that $|a|\le A$ and $|b| \le B$, where
$A$ and $B$ are small relative to~$x$.  More precisely, we
investigate behavior with respect to the Sato--Tate conjecture,
cyclicity, and divisibility of the number of points by a fixed
integer $m$.
\end{abstract}

\section{Introduction}
\label{sec:intro}

\subsection{Motivation}

For integers $a$ and $b$ such that $4a^3+27b^2\ne 0$, we denote by
$\E_{a,b}$ the elliptic curve defined by the \emph{affine
Weierstra\ss\ equation}:
$$
\E_{a,b}~:~Y^2 = X^3 + aX + b.
$$
For a basic background on elliptic curves, we refer the reader to
the book~\cite{Silv} by Silverman.

For any prime $p>3$, we denote by $\Fp$ the finite field with $p$
elements, which we identify with the set of integers $\{0,\pm
1,\ldots,\pm (p-1)/2\}$.

When $p\nmid 4a^3+27b^2$, the set $\E_{a,b}(\Fp)$ consisting of
the $\Fp$-rational points of $\E_{a,b}$ together with a point at
infinity forms an \emph{abelian group} under an appropriate
composition rule called \emph{addition}, and the number of
elements in the group $\E_{a,b}(\Fp)$ satisfies the \emph{Hasse
bound}:
$$
\bigl|\#\E_{a,b}(\F_p)- p - 1\bigr|\le 2\sqrt{p}
$$
(see, for example, \cite[Chapter~V, Theorem~1.1]{Silv}).

A well known conjecture in the theory of elliptic curves, known as
the \emph{Lang--Trotter conjecture} (see~\cite{LaTr1}), asserts
that for any elliptic curve $\E_{a,b}$ and any fixed integer $t$,
the number $\pi_{a,b}(t;x)$ of primes $p\le x$ (with $p\nmid
4a^3+27b^2$) such that
$$
\#\E_{a,b}(\F_p) = p + 1 - t
$$
satisfies the asymptotic formula
$$
\pi_{a,b}(t;x)\sim c_{a,b,t}\cdot\frac{\sqrt x}{\log
x}\qquad(x\to\infty)
$$
with some constant $ c_{a,b,t}$ that depends only on $a$, $b$, and
$t$, provided that $\E_{a,b}$ does not have \emph{complex
multiplication} (see~\cite[Section~III.4]{Silv}) or $t$ is
nonzero. The Lang--Trotter conjecture remains open, although some
progress has been made (see the survey~\cite{MuShp}).

Fouvry and Murty~\cite{FoMu} have studied the problem of
estimating $\pi_{a,b}(0;x)$ \emph{on average} over integers $a,b$
with $|a| \le A$ and $|b| \le B$ and have shown
(see~\cite[Theorem~6]{FoMu}) that the asymptotic formula
$$
\frac{1}{4AB}\sum_{|a|\le A}\sum_{|b|\le
B}\pi_{a,b}(0;x)\sim\frac{\pi}{3}\cdot\frac{\sqrt x}{\log
x}\qquad(x\to\infty)
$$
holds uniformly in the range
\begin{equation}
\label{eq:FM threshold} AB\ge x^{3/2 + \eps} \mand \min\{A,B\} \ge
x^{1/2 + \eps},
\end{equation}
where $\eps>0$ is fixed. For the case $t\ne 0$, David and
Pappalardi~\cite{DavPapp1} have established the following
asymptotic formula in a shorter range of $A$ and $B$:
$$
\frac{1}{4AB}\sum_{|a|\le A}\sum_{|b|\le B}\pi_{a,b}(t;x)\sim
C_t\cdot\frac{\sqrt x}{\log x}\qquad(x\to\infty),
$$
where
$$
C_t=\frac{2}{\pi}\prod_{p\,\mid\,t}\(1-\frac{1}{p^2}\)^{-1}
\prod_{p\,\nmid\,t}\frac{p(p^2-p-1)}{(p-1)(p^2-1)}\,.
$$
This work has been extended by Baier~\cite{Baier1} to the full
range~\eqref{eq:FM threshold}. Other results along these lines
have been obtained
in~\cite{AkDavJur,Baier2,BBIJ,DavPapp2,Gek,James,JamYu}.

We remark that another type of problem ``on average'' over a
similar family of curves has been considered in~\cite{Duke,Grant}.

Here, we investigate the average behavior of the family of curves
$\E_{a,b}$ with $|a|\le A$ and $|b|\le B$ with respect to some
natural statistical properties of their reductions modulo~$p$.
Although these properties are expected to hold for individual
curves, such results remain inaccessible.

\subsection{Our Results}

In the present paper, we study how the family of curves $\E_{a,b}$
with $|a| \le A$ and $|b| \le B$ behaves with respect to:
\begin{itemize}
\item the  \emph{Sato--Tate conjecture} about the distribution of
the cardinalities $\#\E_{a,b}(\F_p)$ (see~\cite{Katz});

\item \emph{cyclicity} of the group $\E_{a,b}(\F_p)$, a notion
which essentially dates back to the work of Borosh, Moreno and
Porta~\cite{BMP} and of Serre~\cite{Ser};

\item \emph{divisibility} of $\#\E_{a,b}(\F_p)$ by a given integer
$m$.
\end{itemize}

Accordingly, for real $0 \le \alpha < \beta \le \pi$, we define
the \emph{Sato--Tate  density}
$$
\mu_{\tt ST}(\alpha,\beta) = \frac{2}{\pi}\int_\alpha^\beta
\sin^2\theta\, d \theta,
$$
and we define the angle $\psi_{a,b}(p) \in [0, \pi]$ via the
identity
$$
p+1 - \#\E_{a,b}(\F_p)  =  2\sqrt{p}\, \cos \psi_{a,b}(p).
$$
We denote by $\Pi^{\tt ST}_{a,b}(\alpha,\beta;x)$ the number of
primes $p \le x$ (with $p \nmid 4a^3+27b^2$) for which $ \alpha
\le \psi_{a,b}(p) \le \beta$. The  \emph{Sato--Tate  conjecture}
asserts that if $\E_{a,b}$ does not have complex multiplication,
then the asymptotic formula
\begin{equation}
\label{eq:ST conj} \Pi^{\tt ST}_{a,b}(\alpha,\beta;x) \sim
\mu_{\tt ST}(\alpha,\beta) \cdot\frac{x}{\log x}\qquad(x\to\infty)
\end{equation}
holds (see~\cite{Birch,Katz,Murt-VK}).

It is well known that $\E_{a,b}(\F_p)$ is an abelian group of rank
at most two. We denote by $\Pi^{\tt C}_{a,b}(x)$ the number of
primes $p \le x$ (with $p \nmid 4a^3+27b^2$) for which
$\E_{a,b}(\F_p)$ is \emph{cyclic}. The conjectured asymptotic
formula
$$
\Pi^{\tt C}_{a,b}(x) \sim C_{a,b}\cdot\frac{x}{\log
x}\qquad(x\to\infty),
$$
where $C_{a,b}$ is a constant that depends only on $a$ and $b$,
has been established conditionally (under the \emph{Extended
Riemann Hypothesis}) in some cases, and there are several
unconditional lower bounds on $\Pi^{\tt C}_{a,b}(x)$; see the
original papers~\cite{Coj1,Coj2,CojMur,CojFouMur,GuMu,Murt-MR,Ser}
as well as the recent surveys~\cite{Coj3,MuShp}.

Finally, for a fixed integer $m\ge 1$, we denote by $\Pi^{\tt
D}_{a,b}(m;x)$ the number of primes $p \le x$ (with $p \nmid
4a^3+27b^2$) for which $m \mid \# \E_{a,b}(\F_p)$.

It worth mentioning that the \emph{Chebotarev density theorem} can
be used to study $\Pi_{a,b}^{\tt D}(m;x)$ for individual curves
(see~\cite{Coj3}). By averaging over $a$ and $b$, however, we
obtain sharper results which are also uniform in $m$ up to any
fixed power of $\log x$.

Taylor~\cite{Tayl} has recently announced a complete proof of the
Sato--Tate conjecture, which implies~\eqref{eq:ST conj} in
particular, but his work on individual curves does not imply any
results on average due to the lack of uniformity with respect to
the coefficients $a$ and $b$ in the Weierstra\ss\ equation.

Here, we obtain an asymptotic formula for the number of pairs
$(a,b)$ with $|a| \le A$, $|b| \le B$ and $p \nmid 4a^3+27b^2$
such that $\E_{a,b}(\F_p)$ belongs to a certain sufficiently
``massive'' collection of isomorphism classes of elliptic curves.
Using this result, we derive asymptotic formulae for the sums
\begin{equation*}
\begin{split}
&N^{\tt ST}_{\alpha, \beta}(A,B;x) = \sum_{|a|\le A}\sum_{|b|\le
B}
\Pi^{\tt ST}_{a,b}(\alpha,\beta;x), \qquad \qquad \\
&N^{\tt C}(A,B;x)=\sum_{|a|\le A}\sum_{|b|\le B}\Pi^{\tt C}_{a,b}(x), \\
&N^{\tt D}_m(A,B;x)=\sum_{|a|\le A}\sum_{|b|\le B}\Pi^{\tt
D}_{a,b}(m;x).
\end{split}
\end{equation*}

The main terms in our expansions of $N^{\tt ST}_{\alpha,
\beta}(A,B;x)$, $N^{\tt C}(A,B;x)$ and $N^{\tt D}_m(A,B;x)$ are
derived from asymptotic formulae of Birch~\cite{Birch},
Vl\u{a}du\c{t}~\cite{Vlad} and Howe~\cite{Howe}, respectively. The
asymptotic formula of Birch~\cite{Birch} is not quite sufficient
for our purposes, however, so we give an explicit bound for the
error term which is obtained using the method of
Niederreiter~\cite{Nied}.

In the case of the Sato--Tate distribution, the computation of the
error term is almost trivial. The other cases require a more
specialized treatment, and these are done using now standard
techniques; for example, we apply a result of Indlekofer, Wehmeier
and Lucht~\cite{IWL}.

In the last section, we give an outline of several other questions
concerning reductions of elliptic curves that can be approached
with our method.

\subsection{Our Method}

The functions $N^{\tt ST}_{\alpha, \beta}(A,B;x)$, $N^{\tt
C}(A,B;x)$ and $N^{\tt D}_m(A,B;x)$ can be studied via the method
of Fouvry and Murty~\cite{FoMu}, which makes essential use of the
\emph{Weil bound}; see~\cite[Chapter~5]{LN}. Here, however, we
obtain sharper results by applying  bounds on multiplicative
character sums rather than estimating exponential sums as
in~\cite{FoMu}. Using the \emph{Polya--Vinogradov} and
\emph{Burgess bounds} (see~\cite[Theorems~12.5 and 12.6]{IwKow})
one already obtains stronger results for individual primes than
with exponential sums. Moreover, the use of multiplicative
character sums allows for additional savings on average over
primes $p\le x$. In the present paper, we use a result of
Garaev~\cite{Gar} on multiplicative character sums, which is
derived from a variant of the large sieve inequality
(see~\cite[Section~7.4]{IwKow}), and we obtain nontrivial bounds
in a region that is significantly wider than~\eqref{eq:FM
threshold}.
Our method leads to bounds which are nontrivial whenever $A$ and
$B$ satisfy the inequalities
\begin{equation}
\label{eq:BS threshold 1} A,B\ge x^\eps\mand AB\ge x^{1+\eps}.
\end{equation}
However, to avoid some tedious technical complications and more
cluttered expressions for the error term, we further assume that
\begin{equation}
\label{eq:BS threshold 2} A,B \le x^{1-\eps}
\end{equation}
(which of course makes the first inequality in~\eqref{eq:BS
threshold 1} redundant). Certainly, the case in which $A$ and $B$
are both small is of primary interest, so the
restriction~\eqref{eq:BS threshold 2} is rather mild. 

One of the main ingredients of the method of Fouvry and Murty is
the use of the Weil bound to prove the asymptotic formula $2AB/p +
O(p^{1/2 + o(1)})$ for the number of curves $\E_{a,b}$ with $|a|
\le A\le (p-1)/2$ and $|b|\le B\le (p-1)/2$ that are isomorphic to
a given curve $\E_{r,s}$; see~\cite[Section~7]{FoMu}. Here, we
show that, on average over $r$ and $s$, the error term can be
improved substantially, and this suffices for the problems that we
consider. On the other hand, our method does not directly apply to
the question considered in~\cite{FoMu} since a set of elliptic
curves over $\F_p$ with a prescribed number of
\hbox{$\F_p$-rational} points (that is, a set of isogenous curves)
is much ``thinner'' than the sets of curves with which we work. Of
course, there is some possibility that both approaches might be
combined to improve the threshold~\eqref{eq:FM threshold} for the
original problem.

Baier and Zhao~\cite{BaZh} have also studied the distribution of
$N^{\tt ST}_{\alpha, \beta}(A,B;x)$ using a very different method
from ours; their results are also different (but there is partial
overlap) and in many cases are weaker with respect to the range of
$A$ and $B$ as well as the uniformity in $\alpha$ and $\beta$. In
particular, among other restrictions, the inequalities
$$
A,B \ge x^{1/2+\eps}
$$
are required for~\cite[Theorem~1]{BaZh}. In some cases, however,
the results of~\cite{BaZh} are stronger than ours.  It is worth
mentioning that Baier and Zhao~\cite{BaZh} have estimated the
average deviation of $\Pi^{\tt ST}_{a,b}(\alpha,\beta;x)$ from the
value predicted by the Sato--Tate conjecture.

Finally, several more results in these directions
  have recently been obtained in~\cite{Baier2}
that also appear to be weaker than our results.

\subsection{Notation}

Throughout the paper, any implied constants in the symbols $O$ and
$\ll$ may occasionally depend, where obvious, on the parameters
$\eps$ and  $K$ but are absolute otherwise. We recall that the
notations $U \ll V$ and  $U = O(V)$ are both equivalent to the
statement that the inequality $|U| \le c\,V$ holds with some
constant $c> 0$.

The letters $p$ and $q$ always denote prime numbers, while $m$ and
$n$ always denote integers. As usual, we use $\pi(x)$ to denote
the number of primes $p\le x$.

\subsection{Acknowledgements}

The authors are grateful to Antal Balog for fruitful discussions
which have led to an improvement of our original results. The
authors would also like to thank Nick Katz for several comments
concerning the Sato--Tate conjecture and in particular for his
suggestion of Lemma~\ref{lem:ST Stat}. This work began during a
pleasant visit by W.~B.\ to Macquarie University; the support and
hospitality of this institution are gratefully acknowledged.
During the preparation of this paper, I.~S.\ was supported in part
by ARC grant DP0556431.

\section{Preliminaries}
\label{sec:prep}

\subsection{Character sums}

For a prime $p$, we denote by $\cX_p$ the set of multiplicative
characters of $\Fp$, $\chi_0$ the principal character of $\Fp$,
and $\cX_p^*=\cX_p\setminus\{\chi_0\}$ the set of nonprincipal
characters; we refer the reader to~\cite[Chapter~3]{IwKow} for the
necessary background on multiplicative characters. We recall the
following orthogonality relations:
\begin{equation}
\label{eqn:orth chi/u} \frac{1}{p-1} \sum_{\chi \in\cX_p}\chi(v) =
\left\{\begin{array}{ll}
1&\quad\text{if $v=1$;}\\
0&\quad\text{otherwise,}
\end{array}\right.
\end{equation}
and
\begin{equation}
\label{eqn:orth chi_1/chi_2} \frac{1}{p-1} \sum_{u\in
\F_p^*}\chi_1(u) \overline\chi_2(u) = \left\{\begin{array}{ll}
1&\quad\text{if $\chi_1=\chi_2$;}\\
0&\quad\text{otherwise,}
\end{array}\right.
\end{equation}
for all $v \in \F_p$ and  $\chi_1,\chi_2 \in\cX_p$ (here,
$\overline\chi_2$ is the character obtained from $\chi_2$ by
complex conjugation).

The following result combines the Polya--Vinogradov bound (for
$\nu=1$) with the Burgess bounds (for $\nu\ge 2$);
see~\cite[Bound~(12.58)]{IwKow} and also~\cite[Theorems~12.5 and
12.6]{IwKow}:

\begin{lemma}
\label{lem:PVB} Uniformly for all primes $p$, all positive
integers $L,M,\nu$, and all characters $\chi\in\cX_p^*$, we have
$$
\sum_{n=L+1}^{L+M} \chi(n)  \ll  M^{1 -1/\nu} p^{(\nu+1)/(4\nu^2)}
(\log p)^{1/\nu}.
$$
\end{lemma}

The next bound is due to Ayyad, Cochrane and
Zheng~\cite[Theorem~2]{ACZ}; see also the result of Friedlander
and Iwaniec~\cite{FrIw}:

\begin{lemma}
\label{lem:Sum 4} Uniformly for all positive integers $L,M$, we
have
$$
\sum_{\chi \in \cX_p^*} \left|\,\sum_{n=L+1}^{L+M} \chi(n)
\right|^4 \ll  p M^{2+ o(1)}\qquad(p\to\infty).
$$
\end{lemma}

We also need the following statement, which is contained in the
more general result~\cite[Theorem~10]{Gar} of Garaev (which also
applies to character sums with composite moduli and allows
significantly more flexibility in the choice of $M$):

\begin{lemma}
\label{lem:Gar} Fix $\varepsilon>0$ and $\eta>0$. If $x$ is
sufficiently large, then for all $M\ge x^\varepsilon$, all primes
$p\le x$ with at most $x^{3/4 + 4\eta +o(1)}$ exceptions as
$x\to\infty$, and all characters $\chi \in \cX_p^*$, we have
$$
\left|\,\sum_{n= 1}^{ M} \chi(n)\right|  \le M^{1 -\eta} ,
$$
where the function implied by $o(1)$ depends only on $\varepsilon$
and $\eta$.
\end{lemma}

\begin{proof} We can assume that
$\eta<1/16$ for otherwise there is nothing to prove.

If $M \ge x^{3/4}$, then the result follows from
Lemma~\ref{lem:PVB} with $\nu=1$. Indeed, for a prime $p\le x$ we
have
$$
\left|\,\sum_{n=1}^{ M} \chi(n)\right| \le p^{1/2+o(1)}   \le
x^{1/2+o(1)}\le M^{2/3 + o(1)} \le M^{1-\eta}
$$
for any fixed $\eta\in(0,1/16)$ if $x$ is large enough.

For $M\le x^{3/4}$ the result is a direct consequence
of~\cite[Theorem~10]{Gar}.
\end{proof}

\subsection{Distribution of powers}

Let $d_p = \gcd(p-1,6)$ and put
\begin{equation}
\label{eq:def sigma} \sigma_p(M) = \max_{\substack{\chi\in
\cX_p^*\\ \chi^{d_p}=\chi_0}}\left\{1,~\left|\,\sum_{n= 1}^{ M}
\chi(n)\right|\right\}.
\end{equation}

For any integers $B,s$ we define
$$
\cZ_s(B;p)=\{u\in\Fp^*~:~s u^6 \equiv b \pmod p   \text{ where }
|b|\le B\}.
$$
We have the following bound on the cardinality of $\cZ_s(B;p)$:

\begin{lemma}
\label{lem:ZsB} For all primes $p$ and all positive integers
$B,s<p$, we have
$$
\bigl|\#\cZ_s(B;p) - 2B \bigr|\le 11\,\sigma_p(B).
$$
\end{lemma}

\begin{proof}
For all $n\in\Z$ we have
$$
\#\{u\in\Fp^*~:~u^6\equiv n\pmod p\}=\sum_{\substack{\chi\in
\cX_p\\ \chi^{d_p}=\chi_0}}\chi(n).
$$
If $\overline s$ is an integer such that $s\overline s\equiv
1\pmod p$, it follows that
$$
\#\cZ_s(B;p)=\sum_{|b|\le B}\sum_{\substack{\chi\in \cX_p\\
\chi^{d_p}=\chi_0}}\chi(\overline s b)
=(2B+1)+\sum_{\substack{\chi\in \cX_p^*\\ \chi^{d_p}=\chi_0}}
\overline \chi(s) \sum_{|b|\le B}\chi(b).
$$
Since the inner sum is bounded by
$$
\left|\,\sum_{|b|\le B} \chi(b)\right|  \le 2\,\sigma_p(B),
$$
and
$$
\#\{\chi \in \cX_p^* :~\chi^{d_p}=\chi_0\}=d_p-1\le 5,
$$
the result follows.
\end{proof}

For any integers $A,B,r,s$ we define
$$
\cZ_{r,s}(A,B;p)=\{u\in \cZ_s(B;p)~:~r u^4 \equiv a \pmod p \text{
where } |a|\le A \}.
$$

\begin{lemma}
\label{lem:ZrsAB-1} For all primes $p$ and all positive integers
$A,B,s<p$, we have
$$
\sum_{r\in\Fp}\left|\#\cZ_{r,s}(A,B;p) - \frac{2 A\cdot\#
\cZ_s(B;p)}{p}\right|\le A^{1/2} B p^{1/4+o(1)} + A^{1/2} B^{1/2}
p^{1/2 + o(1)}
$$
as $p\to\infty$.
\end{lemma}

\begin{proof} We can assume that $AB>p$ since the result is trivial
otherwise.
Indeed, if $AB\le p$ then $A^{1/2} B^{1/2} p^{1/2}\ge AB$, while
\begin{equation*}
\begin{split}
\sum_{r\in\Fp}&\left|\#\cZ_{r,s}(A,B;p) - \frac{2 A\cdot\#
\cZ_s(B;p)  }{p} \right|\\ & \qquad  \le
\sum_{r\in\Fp}\#\cZ_{r,s}(A,B;p) + 2 A\cdot\# \cZ_s(B;p)  \ll AB.
\end{split}
\end{equation*}

For every $a\in\F_p^*$ let $\overline a$ be an integer such that
$a\overline a\equiv 1\pmod p$.  Using~\eqref{eqn:orth chi/u} it
follows that
\begin{equation*}
\begin{split}
\#\cZ_{r,s}&(A,B;p)=\sum_{u\in \cZ_s(B;p) } \sum_{0 < |a|\le A}
\frac{1}{p-1} \sum_{\chi\in\cX_p}\chi(ru^4 \overline{a})\\
&=\frac{2 A\cdot\# \cZ_s(B;p)  }{p-1} +O(1)+
\frac{1}{p-1}\sum_{\chi\in\cX_p^*} \chi(r) \sum_{u\in \cZ_s(B;p)
}\chi(u^4) \sum_{|a|\le A} \overline\chi(a).
\end{split}
\end{equation*}
Since
$$
\frac{2A\cdot\# \cZ_s(B;p)}{p}-\frac{2 A\cdot\#\cZ_s(B;p)}{p-1}
\ll\frac{A\cdot\#\cZ_s(B;p)}{p^2}\ll\frac{AB}{p^2}\ll 1,
$$
we have
\begin{equation}
\label{eq:Mult Sum} \sum_{r\in\Fp}\left|\#\cZ_{r,s}(A,B;p)  -
\frac{2 A\cdot\# \cZ_s(B;p)  }{p} \right| \ll  p +  W,
    \end{equation}
where
$$
W= \frac{1}{p} \sum_{r\in\Fp} \left|\,\sum_{\chi\in\cX_p^*}
\chi(r) \sum_{u\in \cZ_s(B;p) }\chi(u^4) \sum_{|a|\le A}
\overline\chi(a) \right|.
$$
By the Cauchy inequality, $W^2$ does not exceed
\begin{equation*}
\begin{split}
  \frac{1}{p}&\sum_{r\in\Fp}\left|\,\sum_{\chi\in\cX_p^*} \chi(r)
\sum_{u\in \cZ_s(B;p) }\chi(u^4) \sum_{|a|\le A} \overline\chi(a)
\right|^2\\
&=\frac{1}{p} \sum_{\chi_1, \chi_2\in\cX_p^*} \sum_{u_1,u_2\in
\cZ_s(B;p)} \chi_1(u_1^4) \overline\chi_2(u_2^4) \sum_{|a_1|,
|a_2|\le A} \overline\chi_1(a_1) \chi_2(a_2)\sum_{r\in\Fp}
\chi_1(r) \overline\chi_2(r).
\end{split}
\end{equation*}
Using the orthogonality relation~\eqref{eqn:orth chi_1/chi_2} we
deduce that
\begin{equation}
\label{eqn:W2} W^2  \le   \sum_{\chi\in\cX_p^*} \left|\sum_{u\in
\cZ_s(B;p)} \chi(u^4) \right|^2
    \left|\sum_{|a| \le A}  \chi(a)\right|^2.
\end{equation}
Applying the  Cauchy inequality again, it follows that
\begin{equation}
\label{eqn:W4} W^4  \le   \sum_{\chi\in\cX_p^*} \left|\sum_{u\in
\cZ_s(B;p)} \chi(u^4) \right|^4 \cdot\sum_{\chi\in\cX_p^*}
\left|\sum_{|a| \le A}  \chi(a)\right|^4.
\end{equation}
The second sum is of size $O(p^{1+o(1)}A^2)$ by Lemma~\ref{lem:Sum
4}. For the first sum, we extend the summation to include the
trivial character $\chi=\chi_0$, obtaining
$$
\sum_{\chi\in\cX_p^*} \left|\sum_{u\in \cZ_s(B;p)} \chi(u^4)
\right|^4 \le \sum_{\chi\in\cX_p} \left|\sum_{u\in\cZ_s(B;p)}
\chi(u^4) \right|^4 = p\,T,
$$
where $T$ is the number of solutions to the congruence
$$
u_1^4u_2^4 \equiv u_3^4 u_4^4 \pmod p,\qquad u_1,u_2,u_3,u_4\in
\cZ_s(B;p).
$$
Note that $T$ does not exceed the number of quadruples
$(u_1,u_2,u_3,u_4)$ in $\cZ_s(B;p)^4$ for which
$$
u_1^{12}u_2^{12}\equiv u_3^{12}u_4^{12}\pmod p.
$$
Since $su_j^6 \equiv b_j\pmod p$ for some $b_j$ with $|b_j|\le B$,
and each $b_j$ corresponds to at most six values of $u_j$, it
follows that $T\le 6^4 R$, where $R$ is the number of solutions to
the congruence
$$
b_1^2b_2^2 \equiv b_3^2 b_4^2 \pmod p, \qquad |b_1|, |b_2|, |b_3|,
|b_4| \le B.
$$
Clearly, $R = 2\,Q$, where $Q$ is the number of solutions to the
congruence
$$
b_1b_2\equiv b_3b_4\pmod p,\qquad |b_1|,|b_2|,|b_3|,|b_4|\le B.
$$
Writing
$$
Q = \frac{1}{p-1}\sum_{\chi\in\cX_p} \left|\sum_{|b| \le B}
\chi(b) \right|^4 =  \frac{(2B +1)^4}{p-1} +
\frac{1}{p-1}\sum_{\chi\in\cX_p^*} \left|\sum_{|b| \le B} \chi(b)
\right|^4
$$
and using Lemma~\ref{lem:Sum 4} again, we see that
$$
T \ll R \ll Q \ll B^4p^{-1}  + B^2p^{o(1)}.
$$
Collecting the above estimates and substituting them
into~\eqref{eqn:W4} we deduce that
$$
W^4 \ll p^{2+o(1)} A^2 (B^4p^{-1}+B^2)
$$
which together with~\eqref{eq:Mult Sum} implies that
\begin{equation*}
\begin{split}
\sum_{r\in\Fp}&\left|\#\cZ_{r,s}(A,B;p) - \frac{2 A\cdot\#
\cZ_s(B;p)  }{p} \right| \\
&\qquad  \ll p +  A^{1/2} B p^{1/4+o(1)} + A^{1/2} B^{1/2} p^{1/2
+ o(1)} .
\end{split}
\end{equation*}
Finally, for $AB>p$ we have $p<A^{1/2}B^{1/2}p^{1/2}$, and the
result follows.
\end{proof}

Combining Lemmas~\ref{lem:ZsB} and~\ref{lem:ZrsAB-1} we
immediately obtain:

\begin{cor}
\label{cor:ZrsAB-1} For all primes $p$ and all positive integers
$A,B,s < p$, we have
\begin{equation*}
\begin{split}
\sum_{r\in\Fp}& \left|\#\cZ_{r,s}(A,B;p) - \frac{4AB}{p} \right|\\
& \qquad  \ll A\,\sigma_p(B)+A^{1/2} B p^{1/4+o(1)} + A^{1/2}
B^{1/2} p^{1/2 + o(1)}.
\end{split}
\end{equation*}
\end{cor}

For large values of $A$, the following lemma provides a stronger
bound for the sum considered in Lemma~\ref{lem:ZrsAB-1}.

\begin{lemma}
\label{lem:ZrsAB-2} For all primes $p$ and all positive integers
$A,B,s<p$, we have
$$
\sum_{r\in\Fp}\left|\#\cZ_{r,s}(A,B;p)-
\frac{2A\cdot\#\cZ_s(B;p)}{p}\right|\ll B^{1/2} p\,\log p.
$$
\end{lemma}

\begin{proof} We argue as in the proof of  Lemma~\ref{lem:ZrsAB-1}
arriving at~\eqref{eqn:W2}. Next, we apply Lemma~\ref{lem:PVB}
with $\nu=1$ (that is, the Polya--Vinogradov bound) followed by
the Cauchy inequality, deriving the bound
\begin{equation*}
\begin{split}
W^2&\ll p(\log p)^2\sum_{\chi\in\cX_p^*}\left|\sum_{u\in
\cZ_s(B;p)}\chi(u^4)\right|^2\\
&\le p(\log p)^2\sum_{\chi\in\cX_p}\left|\sum_{u\in
\cZ_s(B;p)}\chi(u^4)\right|^2=p(p-1)(\log p)^2T,
\end{split}
\end{equation*}
where $T$ is the number of solutions to the congruence
$$
u_1^4  \equiv u_2^4  \pmod p, \qquad u_1, u_2\in \cZ_s(B;p).
$$
Since $T\le 4\#\cZ_s(B;p)=O(B)$, the result follows.
\end{proof}

We remark that for in the proof of Lemma~\ref{lem:ZrsAB-2} one can
use Lemma~\ref{lem:PVB} with values of $\nu$ other than one, but
doing so does not lead to any improvement over the bound of
Lemma~\ref{lem:ZrsAB-1}.

Combining Lemmas~\ref{lem:ZsB} and~\ref{lem:ZrsAB-2} we obtain:

\begin{cor}
\label{cor:ZrsAB-2} For all primes $p$ and all positive integers
$A,B,s < p$, we have
$$
\sum_{r\in\Fp}\left|\#\cZ_{r,s}(A,B;p) - \frac{4AB}{p} \right| \ll
A\,\sigma_p(B)+ B^{1/2} p\,\log p  .
$$
\end{cor}

\subsection{Statistics of elliptic curves}
\label{sec:Ell Curves}

It is well known that if $a,b,r,s\in\Fp$, then the two curves
$\E_{a,b}$ and $\E_{r,s}$ are \emph{isomorphic over $\F_p$} if and
only if $a=ru^4$ and $b=su^6$ for some $u \in \F_p^*$. In
particular, each curve $\E_{a,b}$ with $a,b\in \F_p^*$ is
isomorphic to $(p-1)/2$ elliptic curves $\E_{r,s}$, and there are
$2p + O(1)$ distinct isomorphism classes of elliptic curves over
$\Fp$; see~\cite{Len}. Thus, our results can be conveniently
formulated in terms  of counting functions for individual curves
$\E_{a,b}$ rather than in terms of isomorphism classes of curves,
as in the papers~\cite{Birch,Howe,Vlad}.

Let $\cT_p(\alpha,\beta)$ be the set of set of pairs $(a,b)\in
\F_p^*\times\F_p^*$ such that  the inequalities $ \alpha \le
\psi_{a,b}(p) \le \beta$ hold. Thanks to Birch~\cite{Birch}, one
knows that
$$
\# \cT_p(\alpha,\beta) \sim \mu_{\tt ST}(\alpha,\beta) p^2
\qquad(p\to\infty),
$$
however we require a stronger result. What is needed is a full
analogue for the Sato--Tate density of the bound of
Niederreiter~\cite{Nied} on the discrepancy in the distribution of
values of (normalized) Kloosterman sums. Fortunately, such a
result can be obtained using the same methods since all of the
underlying tools, namely~\cite[Lemma~3]{Nied}
and~\cite[Theorem~13.5.3]{Katz}, apply to $\psi_{a,b}(p) $ as well
as to values of Kloosterman sums. In particular,
from~\cite[Theorem~13.5.3]{Katz} it follows that
$$
\frac{1}{(q-1)^2} \sum_{\substack{a,b\in \F_p^*\\4a^3+27b^2 \ne
0}} \frac{\sin\((n+1)\psi_{a,b}(p)\)}{\sin\( \psi_{a,b}(p)\)} \ll
nq^{-1/2}  \qquad (n =1,2, \ldots\,)
$$
(see also the work of Fisher~\cite[Section~5]{Fish}). Thus, as
in~\cite{Nied}, we have:

\begin{lemma}
\label{lem:ST Stat} Uniformly for all primes $p$, we have
$$
\max_{0 \le \alpha < \beta \le \pi} \left|\# \cT_p(\alpha,\beta) -
\mu_{\tt ST}(\alpha,\beta) p^2 \right| \ll p^{7/4}.
$$
\end{lemma}

Next, for any prime $p$ we denote
$$
\vartheta_p=\prod_{q\,\mid\,p-1}\(1-\frac{1}{q\(q^2-1\)}\),
$$
where the product is taken over all prime divisors $q$ of $p-1$.

Let $\cC_p$ be the set of pairs $(a,b)\in \F_p^*\times\F_p^*$ such
that $\E_{a,b}(\Fp)$ is cyclic. The cardinality of $\cC_p$ has
been estimated by Vl\u{a}du\c{t}~\cite{Vlad} as follows:

\begin{lemma}
\label{lem:Cycl Stat} For all primes $p$, we have
$$
\left|\# \cC_p - \vartheta_p\,p^2 \right| \le p^{3/2 +
o(1)}\qquad(p\to\infty).
$$
\end{lemma}

Finally, for any integer $k$, let $\omega_k(m)$ denote the
completely multiplicative function which is defined on prime
powers $q^j$ as follows:
\begin{equation}
\label{eq:omega} \omega_k(q^j) = \left\{\begin{array}{ll}
\displaystyle{\frac{1}{q^{j-1}(q-1)} } &\quad\hbox{if $k\not
\equiv 1 \pmod {q^{\rf{j/2}}}$;}\\ \\
          \displaystyle{\frac{q^{\fl{j/2}+1} + q^{\fl{j/2}} -1 }{q^{j +
\fl{j/2}-1}(q^2-1)}}
&\quad\hbox{if $k \equiv 1 \pmod {q^{\rf{j/2}}}$.}\\
\end{array}\right.
\end{equation}

For each integer $m$, let $\cD_p(m)$ be the set of pairs $(a,b)\in
\F_p^*\times\F_p^*$ such that $m \mid \# \E_{a,b}(\Fp)$. Then, by
the result of Howe~\cite{Howe} we have the following asymptotic
formula for $\# \cD_p(m)$:

\begin{lemma}
\label{lem:Div Stat} For all primes $p$ and all positive integers
$m$, we have
$$
\left|\# \cD_p(m) - \omega_p(m) p^2 \right| \le m^{1 + o(1)}
p^{3/2}\qquad(m\to\infty).
$$
\end{lemma}

\section{Main Results}
\label{sec:main}

\subsection{Distribution of curves over finite fields}

For an arbitrary subset $\cS \subseteq \Fp\times\Fp$, we denote by
$M_p(\cS,A,B)$ the number of curves $\E_{a,b}$ such that $(a,b)
\in \cS$, $|a| \le A$ and $|b| \le B$. Here, we obtain an
asymptotic formula for $M_p(\cS,A,B)$.

Similar to~\eqref{eq:def sigma}, we now define $e_p = \gcd(p-1,4)$
and put
$$
\rho_p(M) = \max_{\substack{\chi\in \cX_p^*\\ \chi^{e_p}=\chi_0}}
\left\{1,~\left|\,\sum_{n=1}^{M} \chi(n)\right|\right\}.
$$
We also denote
\begin{equation*}
\begin{split}
\cE_1(A,B;p)&=\min\bigl\{ A\,\sigma_p(B)+   A^{1/2} B p^{1/4}+
A^{1/2} B^{1/2} p^{1/2}, \\
& \qquad \qquad \qquad \qquad B\,\rho_p(A) + A B^{1/2} p^{1/4} +
A^{1/2} B^{1/2} p^{1/2}\bigr\},\\
\cE_2(A,B;p)&=\min\left\{ A\,\sigma_p(B)+ B^{1/2} p\,\log
p,~B\,\rho_p(A) + A^{1/2}p\,\log p\right\}.
\end{split}
\end{equation*}

\begin{theorem}
\label{thm:Set S}  For all primes $p>3$, all integers $1\le A,B\le
(p-1)/2$, and all subsets $\cS \subseteq \Fp\times\Fp$ such that
whenever $(r,s)\in\cS$ and $\E_{a,b}(\F_p)\cong\E_{r,s}(\F_p)$ it
follows that $(a,b)\in\cS$, the following bound holds uniformly:
$$
M_p(\cS,A,B) - \frac{4 AB}{p^2}\,\#\cS  \ll
\min\left\{\cE_1(A,B;p)p^{o(1)},~\cE_2(A,B;p)\right\}.
$$
\end{theorem}

\begin{proof}
It follows from the properties of isomorphic curves given in
Section~\ref{sec:Ell Curves} that
$$
M_p(\cS,A,B)= \frac{1}{p-1} \sum_{(r,s)\in
\cS}\#\cZ_{r,s}(A,B;p)+O(A+B),
$$
where we have estimated the contribution from curves with $ab=0$
trivially as $O(A+B)$; note that if $a=ru^4$ and $b=su^6$ then the
same relations also hold with $-u$ instead of $u$, so each group
$\E_{a,b}(\Fp)$ with $|a|\le A$ and $|b|\le B$ is counted
precisely $p-1$ times in the sum on the right-hand side. Applying
Corollary~\ref{cor:ZrsAB-1}, we obtain that
\begin{equation*}
\begin{split}
M_p(\cS,A,B)= \frac{4AB}{p(p-1)} & \,\#\cS \\
      +~& O\(A \,\sigma_p(B) +
A^{1/2} B p^{1/4+o(1)}+ A^{1/2} B^{1/2} p^{1/2+o(1)}\).
\end{split}
\end{equation*}
Since
$$
\frac{4AB}{p(p-1)}\,\#\cS-\frac{4AB}{p^2}\,\#\cS\ll\frac{AB}{p}
\le A\le A\,\sigma_p(B),
$$
it follows that
$$
M_p(\cS,A,B) - \frac{4 AB}{p^2}\,\#\cS\ll A \,\sigma_p(B) +
A^{1/2} B p^{1/4+o(1)}+ A^{1/2} B^{1/2} p^{1/2+o(1)}.
$$
Examining our arguments closely, in particular those of
Section~\ref{sec:prep}, we see that the roles of $A$ and $B$ are
fully interchangeable, hence we also have
$$
M_p(\cS,A,B) - \frac{4 AB}{p^2}\,\#\cS  \ll B\,\rho_p(A)+AB^{1/2}
p^{1/4+o(1)}+ A^{1/2} B^{1/2} p^{1/2+o(1)},
$$
and thus,
$$
M_p(\cS,A,B)-\frac{4 AB}{p^2}\,\#\cS \ll \cE_1(A,B;p)p^{o(1)}.
$$
Using Corollary~\ref{cor:ZrsAB-2} instead of
Corollary~\ref{cor:ZrsAB-1}, a similar argument shows that
$$
M_p(\cS,A,B)-\frac{4 AB}{p^2}\,\#\cS \ll \cE_2(A,B;p),
$$
and this concludes the proof.
\end{proof}

Using  Lemma~\ref{lem:PVB} to estimate $\cE_1(A,B;p)$,  we deduce
that:

\begin{cor}
\label{cor:Set S Burgess E1} Under the hypotheses of
Theorem~\ref{thm:Set S}, the bound
\begin{equation*}
\begin{split}
&\left|M_p(\cS,A,B) - \frac{4 AB}{p^2}\,\#\cS\right|\\
&\qquad \le \min\Bigl\{ A B^{1 -1/\nu} p^{(\nu+1)/(4\nu^2)} +
A^{1/2}Bp^{1/4}+A^{1/2} B^{1/2} p^{1/2 }, \\
& \qquad \qquad \qquad \qquad A^{1 -1/\nu} B p^{(\nu+1)/(4\nu^2)}
+ AB^{1/2} p^{1/4} + A^{1/2} B^{1/2} p^{1/2 }\Bigr\} p^{o(1)}
\end{split}
\end{equation*}
holds with any fixed integer $\nu\ge 1$, where the function
implied by $o(1)$ depends only on $\nu$.
\end{cor}

Similarly, using  Lemma~\ref{lem:PVB} to estimate $\cE_2(A,B;p)$
we obtain:

\begin{cor}
\label{cor:Set S Burgess} Under the hypotheses of
Theorem~\ref{thm:Set S}, the bound
\begin{equation*}
\begin{split}
&\left|M_p(\cS,A,B) - \frac{4 AB}{p^2}\,\#\cS\right|\\
&\qquad \le \min\left\{ A B^{1 -1/\nu} p^{(\nu+1)/(4\nu^2)} +
B^{1/2} p, ~A^{1 -1/\nu} B p^{(\nu+1)/(4\nu^2)} + A^{1/2}p\right\}
p^{o(1)}
\end{split}
\end{equation*}
holds with any fixed integer $\nu\ge 1$, where the function
implied by $o(1)$ depends only on $\nu$.
\end{cor}

Applying Corollary~\ref{cor:Set S Burgess E1} with a sufficiently
large integer $\nu$, we have the following:

\begin{cor}
Under the hypotheses of Theorem~\ref{thm:Set S}, for any fixed
$\eps>0$ there exists $\delta >0$ such that if
$$
\min\{A,B\} \ge p^{1/4+\eps} \mand AB \ge p^{1+\eps},
$$
then
$$
M_p(\cS,A,B) - \frac{4 AB}{p^2}\,\#\cS  \ll ABp^{-\delta},
$$
where the constant implied by $\ll$ depends only on   $\eps$.
\end{cor}

\subsection{Sato--Tate conjecture on average}

\begin{theorem}
\label{thm:NST}  For any fixed $\eps>0$ there exists $\delta >0$
such that for all integers  $A$ and $B$ satisfying the
inequalities~\eqref{eq:BS threshold 1} and~\eqref{eq:BS threshold
2}, and all real numbers $0 \le \alpha < \beta \le \pi$, we have
$$
N^{\tt ST}_{\alpha, \beta}(A, B;x) = \(4  \mu_{\tt ST}(\alpha,
\beta) + O\bigl(x^{-\delta}\bigr)\)AB\,\pi(x),
$$
where the constant implied by $O$ depends only on  $\eps$.
\end{theorem}

\begin{proof}
Let us assume that $A\ge B$ since the case $A<B$ is similar. Using
the trivial bound $M_p(\cT_p(\alpha, \beta),A,B)\le AB$  for
primes $p\le 2A+1$, we have
$$
N^{\tt ST}_{\alpha, \beta}(A, B;x) =\sum_{2A+1<p\le
x}M_p(\cT_p(\alpha, \beta),A,B)+O(A^2B).
$$
Applying Lemma~\ref{lem:ST Stat} and Theorem~\ref{thm:Set S}, we
derive that
$$
M_p(\cT_p(\alpha, \beta),A,B) - 4\mu_{\tt ST}(\alpha, \beta) AB
\ll \cE_1(A,B;p)p^{o(1)}  + ABp^{-1/4}.
$$
Therefore, since $\mu_{\tt ST}(\alpha, \beta) \ll 1$, we have
$$
    N^{\tt ST}_{\alpha, \beta}(A, B;x)-4\mu_{\tt ST}(\alpha,
\beta)AB\,\pi(x) \ll A^2B + ABx^{3/4}+ x^{o(1)}\sum_{p\le
x}\cE_1(A,B;p).
$$
Applying  Lemma~\ref{lem:Gar} with $\eta = 1/24$,   we get
$$
\sum_{p\le x}\sigma_p(B)\le Bx^{11/12 +o(1)} +B^{23/24}x
\qquad(x\to\infty),
$$
and it follows that
\begin{equation}
\label{eq:aver est}
\begin{split}
\sum_{p\le x}\cE_1(A,B;p)&\ll ABx^{11/12+o(1)}
+AB^{23/24}x\\
&\qquad+A^{1/2} B x^{5/4}+ A^{1/2} B^{1/2}
x^{3/2}\qquad(x\to\infty).
\end{split}
\end{equation}
After simple calculations, we obtain the stated result.
\end{proof}

\subsection{Cyclicity on average}

Let $\varTheta$ denote the following constant:
$$
\varTheta = \prod_{q} \(1 - \frac{1}{q(q -1)\(q^2-1\)}\),
$$
where the product runs over all primes $q$.

\begin{theorem}
\label{thm:NC} Let $\eps>0$ and $K>0$ be fixed. Then, for all
integers $A$ and $B$ satisfying the inequalities~\eqref{eq:BS
threshold 1} and~\eqref{eq:BS threshold 2}, we have
$$
N^{\tt C}(A,B;x) =\(4 \varTheta  +O\( (\log x)^{-K}\)\)
AB\,\pi(x),
$$
where the constant implied by $O$ depends only on   $\eps$ and
$K$.
\end{theorem}

\begin{proof}
Let us assume that $A\ge B$ since the case $A<B$ is similar. Using
the trivial bound $M_p(\cC_p,A,B)\le AB$ for primes $p\le 2A+1$,
we have
$$
N^{\tt C}(A,B;x) =\sum_{2A+1<p\le x}M_p(\cC_p,A,B)+O(A^2B).
$$
Applying Lemma~\ref{lem:Cycl Stat} and Theorem~\ref{thm:Set S}, we
derive that
$$
\bigl|M_p(\cC_p,A,B) - 4\vartheta_p AB\bigr| \le \(\cE_1(A,B;p) +
ABp^{-1/2}\)p^{o(1)}.
$$
Since the bound $AB\le p^2$ implies $ABp^{-1/2}\le
A^{1/2}B^{1/2}p^{1/2}$, and the second term can be dropped; thus,
$$
\bigl|M_p(\cC_p,A,B)-4\vartheta_p AB\bigr|\le\cE_1(A,B;p)p^{o(1)}.
$$
Hence, using~\eqref{eq:aver est}, the inequality $A^2B \le
ABx/(\log x)^{K+1}$, and the trivial bound $\vartheta_p\ll 1$ for
primes $p\le 2A+1$, after simple calculations we derive the
estimate
\begin{equation}
\label{eq:Cest1} N^{\tt C}(A,B;x) =4AB\sum_{p\le x}\vartheta_p
+O\(\frac{ABx}{(\log x)^{K+1}} \).
\end{equation}

Now write
$$
\sum_{p\le x} \vartheta_p = \sum_{p\le x} f(p-1),
$$
where
$$
f(n) = \prod_{q\,\mid\,n} \(1 - \frac{1}{q\(q^2-1\)}\),
$$
the product being taken over all prime divisors $q$ of $n$; note
that $f(n)$ is a multiplicative function. Let $g(n)$ be the
multiplicative function that is defined on prime powers~$q^k$ as
follows:
$$
g(q^k) =  \left\{\begin{array}{ll}
\displaystyle{\frac{-1}{q\(q^2-1\)}}&\quad\hbox{if $k=1$;}\\
0&\quad\hbox{if $k\ge 2$.}
\end{array}\right.
$$
Then,
$$
f(n) = \sum_{d\,\mid\,n}g(d).
$$
It is easy to check that the functions $f(n)$ and $g(n)$ satisfy
the conditions of~\cite[Theorem~3]{IWL}, hence it follows that
\begin{equation}
\label{eq:sumthest} \sum_{p\le x} \vartheta_p = \sum_{p\le x}
f(p-1) = \varTheta\,\pi(x) + O\(\frac{x}{(\log x)^{K+1}}\)
\end{equation}
with
$$
\varTheta=\sum_{n=1}^\infty \frac{g(d)}{\varphi(d)} =\prod_{q}\(1-
\frac{1}{q(q -1)\(q^2-1\)}\).
$$
Inserting the estimate~\eqref{eq:sumthest} into~\eqref{eq:Cest1},
we finish the proof.
\end{proof}

\subsection{Divisibility on average}

Put
\begin{equation}
\label{eq:mu} \mu = \prod_{q^j\,\|\,m} q^{\rf{j/2}},
\end{equation}
and set
$$
\Omega_m = \frac{1}{\varphi(\mu)} \sum_{\substack{1\le k\le
\mu\\\gcd(k,\mu) = 1}}\omega_k(m),
$$
where $\varphi(\mu)$ is the Euler function, and $\omega_k(m)$ is
the completely multiplicative function which is defined on prime
powers $q^j$ by~\eqref{eq:omega}. For example, if $m=q$ is prime,
then we have
$$
\Omega_q = \frac{1}{q-1} \( \frac{q}{q^2-1} + \sum_{k=2}^{q-1}
\frac{1}{q-1}\) = \frac{q^2-2}{(q-1)(q^2-1)}\,.
$$

\begin{theorem}
\label{thm:ND} Let $\eps>0$ and $K>0$ be fixed. Then, for all
integers $A$ and $B$ satisfying the inequalities~\eqref{eq:BS
threshold 1} and~\eqref{eq:BS threshold 2}, and all integers $m
\le (\log x)^K$, we have
$$
N^{\tt D}_{m}(A, B;x) = \(4\,\Omega_m  +O\( (\log x)^{-K}\)\)
AB\,\pi(x) ,
$$
where the constant implied by $O$ depends only on $K$ and $\eps$.
\end{theorem}

\begin{proof}
Let us assume that $A\ge B$ since the case $A<B$ is similar. Using
the trivial bound $M_p(\cD_p(m),A,B)\le AB$ for primes $p\le
2A+1$, we have
$$
N^{\tt D}_{m}(A, B;x)=\sum_{2A+1<p\le x}M_p(\cD_p(m),A,B)+O(A^2B).
$$
Applying Lemma~\ref{lem:Div Stat} and Theorem~\ref{thm:Set S}, we
see that for $p\le x$ and $x\to\infty$:
\begin{eqnarray*}
\bigl|M_p(\cD_p(m),A,B) - 4\,\omega_p(m)AB\bigr|
& \le & \(\cE_1(A,B;p)  + ABp^{-1/2} \)p^{ o(1)} \\
& \le & \cE_1(A,B;p) p^{o(1)},
\end{eqnarray*}
where the second inequality follows from the fact that the term
$ABp^{-1/2}$ never dominates $\cE_1(A,B;p)$ (see the proof of
Theorem~\ref{thm:NC}). Hence, using~\eqref{eq:aver est} and the
trivial bound $\omega_p(m)\ll 1$ for primes $p\le 2A+1$, we
conclude that
\begin{equation*}
\begin{split}
N^{\tt D}_{m}(A, B;x)&=  4AB\sum_{p\le x}\omega_p(m)\\
&\qquad+O\bigl(ABx^{11/12+o(1)} +AB^{23/24}x\\
& \qquad \qquad + A^{1/2} B x^{5/4+o(1)}+ A^{1/2} B^{1/2}
x^{3/2+o(1)} + A^2B\bigr).
\end{split}
\end{equation*}
As the value of $\omega_p(m)$ depends only on the residue class of
$p$ modulo $\mu$, where $\mu$ is given by~\eqref{eq:mu}, using the
\emph{Siegel--Walfisz theorem} (see~\cite[Corollary~5.29]{IwKow})
we immediately obtain the desired result.
\end{proof}

\section{Primes $p$ with $\gcd(p-1,12)=2$}

Clearly, Theorem~\ref{thm:Set S} leads to nontrivial results only
under the condition $\max\{A,B\}\ge p^{1/4+\eps}$ which is
determined by Lemma~\ref{lem:PVB} (this is where the relevant
character sums admit nontrivial estimates). However in the special
case that $\gcd(p-1,12) =2$, using the result of~\cite{BGHBS} on
the density of quadratic residues and nonresidues in short
intervals, one can obtain some nontrivial estimates over a wider
range.

Indeed, if  $\gcd(p-1,6)=2$, then
\begin{eqnarray*}
\cZ_s(B;p)&=&\{u\in\Fp^*~:~s u^6 \equiv b \pmod p   \text{ where }
|b|\le B\}\\
&=&\{u\in\Fp^*~:~s u^2 \equiv b \pmod p   \text{ where } |b|\le
B\}.
\end{eqnarray*}
Thus $\cZ_s(B;p)/2$ is the number of quadratic residues or
nonresidues in the interval $|b|\le B$ (according to whether $s$
is a quadratic residue or nonresidue). Therefore, by the result
of~\cite{BGHBS} we have
$$
\cZ_s(B;p) \gg B
$$
whenever $B \ge p^{1/(4\sqrt{e}\,)+\eps}$ for some fixed
$\varepsilon$ and all sufficiently large $p$.  This new bound can
be used in Lemmas~\ref{lem:ZrsAB-1} and~\ref{lem:ZrsAB-2} as
before. If $\gcd(p-1,4)=2$ as well, then similar arguments can be
applied with respect to $A$. Arguing as in the proof of
Theorem~\ref{thm:Set S}, we obtain a lower bound on $M_p(\cS,A,B)$
which holds under the condition
$$
\min\{A,B\} \ge p^{1/(4\sqrt{e}\,)+\eps} \mand AB \ge p^{1+\eps}.
$$

\section{Further Applications}

For specific ranges of the parameters $A$ and $B$, one can use
Lemma~\ref{lem:PVB} instead of (or in conjunction with)
Lemma~\ref{lem:Gar} to obtain stronger and more explicit bounds
for the error term in Theorem~\ref{thm:NST}. On the other hand, in
Theorems~\ref{thm:NC} and~\ref{thm:ND} the main contribution to
the error comes from the imprecision involved in estimating sums
with $\vartheta_p$ and $\omega_p(m)$, respectively.

Using Lemma~\ref{lem:PVB} in place of Lemma~\ref{lem:Gar} also
allows one to study averages in which the parameters $a$ and $b$
vary over the shifted intervals $[H+1,H+K]$ and $[L+1,L+M]$,
respectively.

Our arguments can also be used to improve the bound
of~\cite{BanShp} on the size of the ``smallest'' Weierstra\ss\
equation which is isomorphic to a given curve $\E$. Given an
elliptic curve $\E$ over  $\F_p$ let us define
$$\mu(\E)=\min\left\{\max\{a,b\}~:~1 \le a,b<p,\ \E_{a,b}\cong \E\right\},
$$
(that is, the minimum is taken over all curves $\E_{a,b}$ which
are isomorphic to $\E$). This question has been recently
considered in~\cite{CQS}, where, using  a variant of the method
of~\cite{FoMu} the bound $\mu(E)=O(p^{3/4})$ has been derived. It
has been shown in~\cite{BanShp} that for almost all curves one has
$\mu(\E) \le p^{2/3 + o(1)}$. Using a variant of
Theorem~\ref{thm:Set S} for the intervals $1 \le a \le A$ and
$1\le b \le B$ one easily derives $\mu(\E)  \le p^{1/2 + o(1)}$
for almost all curves $\E$ over $\F_p$.

Here, we have not used the full strength of the results of
Garaev~\cite{Gar}.  Doing so, one can actually replace the lower
bound $A,B \ge x^{\eps}$ in~\eqref{eq:BS threshold 1} with the
bound $A,B \ge \exp\(c \sqrt{\log x}\,\)$ for an appropriate
constant $c>0$ and drop the condition~\eqref{eq:BS threshold 2}.

For fixed integers $m > k \ge 0$, one can also study the counting
functions $\varpi_{a,b}^{(E)}(m,k;x)$ and
$\varpi_{a,b}^{(t)}(m,k;x)$ of primes $p\le x$ (with $p\nmid
4a^3+27b^2$) such that
$$
\#\E_{a,b}(\F_p) \equiv k \pmod m \mand p+1 - \#\E_{a,b}(\F_p)
\equiv k \pmod m,
$$
respectively. Our method can be adapted to obtain asymptotic
formulae for the average values
$$
\frac{m}{8AB}\sum_{|a|\le A}\sum_{|b|\le B}
\varpi_{a,b}^{(E)}(m,k;x) \mand \frac{m}{8AB}\sum_{|a|\le
A}\sum_{|b|\le B} \varpi_{a,b}^{(t)}(m,k;x)
$$
over a wide range of values in the parameters $A$, $B$ and $m$.

In principle, one can combine our approach with the results
of~\cite{LuSh} to study the distribution of the discriminants of
complex multiplication fields of the curves $\E_{a,b}(\F_p)$, on
average over $a$ and $b$. Such discriminants are related to the
size of the \emph{Tate--Shafarevich group} of $\E_{a,b}(\F_p)$;
thus, it is likely that our approach can be used to improve some
of the estimates of~\cite{CojDuke} on average.

We remark that the distribution of prime values of
$\#\E_{a,b}(\F_p)$ is of great interest in the theory of
cryptography. Our method can be adapted to study this question as
$a$ and $b$ vary over short intervals; see~\cite{BaCoDa}, where
also the challenging task of evaluating the main term has been
treated.

For the problems mentioned above, the corresponding sets of curves
are ``massive'' enough to permit an application of
Theorem~\ref{thm:Set S}; nevertheless, the main obstacle in each
case is the evaluation of the main term, which may require a
significant effort even if the work is feasible.

It is natural to try to combine our approach with recent work of
James and Yu~\cite{JamYu} which studies, on average over $|a| \le
A$ and $|b| \le B$, the number of primes $p\le x$ for which
$p+1-\E_{a,b}(\F_p)$ is a perfect $k$-th power.  In some cases, it
may be possible to lower the threshold on $A$ and $B$. For $k\ge
3$ the corresponding set of curves appears to be too sparse, but
perhaps for $k=2$ there is a chance for our method to yield an
improvement.

\end{document}